# Solving Obstacle Problems using Optimal Homotopy Asymptotic Method

Muhammad Amjad[1], Haider Ali[1*]

1. Department of Mathematics, COMSATS University Islamabad, Vehari Campus, Vehari, 61100, Pakistan.

*Correspondence: haiderali@cuivehari.edu.pk

## Abstract

Differential equations have void applications in several practical situations, sciences, and non-sciences as Euler–Lagrange equation in classical mechanics, Radioactive decay in nuclear physics, Navier–Stokes equations in fluid dynamics, Verhulst equation in biological population growth, Hodgkin–Huxley model in neural action potentials, etc. The cantilever bridge problem is very important in Bridge Engineering and this can be modeled as a homogeneous obstacle problem in Mathematics. Due to this and various other applications, obstacle problems become an important part of our literature. A lot of work is dedicated to the solution of the obstacle problems. However, obstacle problems are not solved by the considered method in the literature we have visited. In this work, we have investigated the finding of the exact solution to several obstacle problems using the optimal homotopy asymptotic method (OHAM). The graphical representation of results represents the symmetry among them.

*Keywords:* Obstacle problems, Optimal homotopy asymptotic method, Numerical problems.

## 1. Introduction

The contemporary height of science and engineering is accomplished because of the development of several applications of Mathematics in modeling physical situations. To remain viable in science and technology in different eras mathematics is used as a basic tool to model the situation, generalize it and develop different techniques to solve these or other models. The invention of high-speed digital computers leads us to formulate problems that can be solved using the techniques of ordinary differential equations (ODEs) and partial differential equations (PDEs).

Variational Inequality theory is the basic branch of mathematics that leads us to model a physical problem as a problem of the differential equation (David Kinderlehrer etc. book). The scheme introduced by the variational inequality theory when used in diverse domains of engineering and math sciences i.e. in the theory of elasticity, solid and liquids mechanics, fluid flow through porous medium, transportation problems, optimal control of machine or production and structural analysis in different sciences, among others,[1–7]. Any imperative establishment for the use of variational diversities is made available using unilateral approaches problem between the elastic bodies and rigid obstacles. According to Kikuchi and Oden [6], the equilibrium model of some given elastic body in contact to a stiff base is deliberated in light of the notions of variation inequality theory (VIT). A large part of the literature is available on the boundary value problems (BVPs) associated with unilateral, obstacle, and contact problems. Our focus in this text is the solution to Obstacle Problems of contact type. We have to know some basic results linked to the phrasing of the obstacle problem.

Let $\mathbb{H}$ be a real Hilbert Space. The inner product and the associated norm are denoted as $\langle .\,,.\rangle$ and $\|.\|$. Let us denote:

$\mathfrak{K}$      Closed Convex set in $\mathbb{H}$

| $\mathfrak{f}$ | Linear continuous functional on $H$ |
|---|---|
| $\mathfrak{T}, \mathfrak{g}$ | Non-Linear Operators from $H \to H$. |

Our problem is to compute $\mathfrak{u} \in \mathfrak{H}$ so that $\mathfrak{g}(\mathfrak{u}) \in \mathfrak{K}$ also satisfies

$$\langle \mathfrak{T}\mathfrak{u}, \mathfrak{g}(\mathfrak{v}) - \mathfrak{g}(\mathfrak{u})\rangle \geq \langle \mathfrak{f}, \mathfrak{g}(\mathfrak{v}) - \mathfrak{g}(\mathfrak{u})\rangle \tag{1'}$$

The inequity (1') is called as General Variational Inequality given in [7]. Solution of (1') exists and is well known in the literature [7 Theorem 3.1, 8, 9 Chapter IV Section 2]. Using (1') we will find the solution to the different obstacle problems in the next section which consider the obstacle problem.

Numerous order obstacles, partial and associated problems can be represented according to the ordinary framework like:

$$\mathfrak{u}^{(\mathfrak{n})}(\mathfrak{x}) = \begin{cases} \mathfrak{P}(\mathfrak{f}(\mathfrak{x}), \mathfrak{g}(\mathfrak{x}), \mathfrak{u}(\mathfrak{x}), \mathfrak{r}), \mathfrak{a} \leq \mathfrak{x} < \mathfrak{c} \\ \mathfrak{Q}(\mathfrak{f}(\mathfrak{x}), \mathfrak{g}(\mathfrak{x}), \mathfrak{u}(\mathfrak{x}), \mathfrak{r}), \mathfrak{c} \leq \mathfrak{x} < \mathfrak{d} \\ \mathfrak{R}(\mathfrak{f}(\mathfrak{x}), \mathfrak{g}(\mathfrak{x}), \mathfrak{u}(\mathfrak{x}), \mathfrak{r}), \mathfrak{d} \leq \mathfrak{x} < \mathfrak{b} \end{cases} \tag{1}$$

by these related boundary values based upon the terms of $n$ determined as

$$\mathfrak{u}^{(\mathfrak{f})}(\mathfrak{a}) = \mathfrak{u}^{(\mathfrak{f})}(\mathfrak{b}) = \alpha_{\mathfrak{f}}, \mathfrak{f} = 0,1,2,\ldots,\mathfrak{n} - 1 \tag{2}$$

$$\mathfrak{u}^{(\mathfrak{f})}(\mathfrak{c}) = \mathfrak{u}^{(\mathfrak{f})}(\mathfrak{d}) = \gamma_{\mathfrak{f}}, \mathfrak{f} = 0,1,2,\ldots,\mathfrak{n} - 1.$$

The expressions $\mathfrak{P}, \mathfrak{Q}$ and $\mathfrak{R}$ are from $\mathbb{R}^4 \to \mathbb{R}$. It is notable that the expression $\mathfrak{u}^{(\mathfrak{f})}(\mathfrak{x}), \mathfrak{f} = 0,1,2,\ldots,\mathfrak{n} - 1$ are continuous on $\mathfrak{c}$ and $\mathfrak{d}$. The variables $\alpha_{\mathfrak{f}}, \beta_{\mathfrak{f}}$ and $\gamma_{\mathfrak{f}}, \mathfrak{f} = 0,1,2,\ldots,\mathfrak{n} - 1$ are true constants (certain definite values). In general, it is hard to compute an analytic structure of the solution of (1) for any option of $\mathfrak{g}(\mathfrak{x})$ and $\mathfrak{f}(\mathfrak{x})$. To this end, a few numeric techniques have been adopted to obtain estimatted results for the problems of this kind (1). This kind of differential equations occurs within impediment. Some of them are associated with partial as well as few with contact problems. We also have many uses in various flashes in applied and pure sciences.

Obstacle problems have a swift enlargement in applications of human life. The last few decades were very important in the construction and refinement of the obstacle problems and their solutions [1, 3]. Various mathematical techniques were developed and used to solve these problems and focused to obtain as much accuracy as we can. Somehow, as in the general form of obstacle problems, there are a lot of situations in which we are unable to find the exact solution to the problem. For solving this, we used numerical techniques and elaborate reasons about the accuracy of the obtained solution. An assortment of numerical methods was formulated and use to discover the mathematical results to the complicated equations together with definite variance methods as well as the problems based on spline. However, such techniques are not only used to solve or designed for solving the complicated questions. In case the included complicated function is recognized by the researcher, hence we will pertain this procedure of Stampacchia and Lewy [7] for distinguish this complicated question through the progression to limit-value questions devoid of limitations utilizing the VIT and penalty function.

The obstacle problem of second-order was discussed and determined in 2001 through the method of cubical spline in [1]. The alike equation of second-order has additionally been resolved via particular parametric cubic spline technique in 2003 [2] and the sextic spline function to develop a numerical technique for solving the system of 2nd order BVPs linked with impediment problems. This work shows that the numerical results computed by the this method are better than those obtained by other methods like finite difference method (FDM), collocation, and the spline methods. A mathematical example was presented to demonstrate the feasible utility of our method [4]. Cubic Lagrange polynomials which associate spline functions to approximate the solutions and error rate is far batter than the existing FDM and spline approach accordingly given in [3] as well as Galerkin's finite element method (FEM) [4] in 2010. The work [5] used the method of B-spline to figure out the similar technique for 2$^{nd}$ order BVPs. In 2013 Adoptive FEM method was applied [12]. The general second-order obstacle problem can be given as:

$$u^{(2)}(x) = \begin{cases} g(x), a \leq x < c \\ f(x)u(x) + g(x) + r, c \leq x < d \\ g(x), d \leq x < b \end{cases} \tag{3}$$

with the associated boundary conditions (BCs) given as $u(a) = u(b) = \alpha_f$. It is given that the functions $u^{(1)}(x)$ is continuous on $c$ and $d$.

Third-order Obstacle Problem was started formally in 1984 in the work of M. Sakai and R.A. Usmani. Quadratic spline results for two-point factors boundary equations related to 3$^{rd}$ order differential problems. Later on, using the penalty function technique, there are a lot of works that find the solution of the 3$^{rd}$ order problems of this kind. This procedure has been used for solving 3$^{rd}$ order obstacle problems by quintic B-spline, FDM, and quartic spline methods in the works [13, 14, 15, 16, 17, 18, 19, 20, 21]. These works cover the era from 1994 to 2006. However, the most generalized work of this case was discussed by A. K. Khalifa and M. A. Noor in "*A numerical approach for odd-order obstacle problems*" in 1994 [19]. The general second-order obstacle problem can be given as:

$$u^{(3)}(x) = \begin{cases} g(x), a \leq x < c \\ f(x)u(x) + g(x) + r, c \leq x < d \\ g(x), d \leq x < b \end{cases} \tag{4}$$

with the associated BCs given as $u(a) = u(b) = \alpha_f$ and $u^{(1)}(a) = u^{(1)}(b) = \beta_f$. It is given that the functions $u^{(1)}(x), u^{(2)}(x)$ are continuous on $b$ and $c$.

In latest years, a few finite distinction and quintic spline strategies were evolved via way of means of [1, 2, 8, 9] to fix the 4$^{th}$ order structures of differential equations related to impediment also PDEs. Khalifa and Noor [5,10] mentioned the opportunity of the use of the collocation technique with quintic splines as primary features for impediment equations. Khan et al. [11] examined a gadget of 4$^{th}$ order impediment BVPs and carried out the parametric quintic spline technique to locate the mathmatical result. Siddiqi and Akram [12–14] advanced the special spline strategies for the answer of the structures of 4$^{th}$ order impediment BVP. Usmani [15] assumed the bending equation of a stretchted orthogonal plate considered across the complete floor through an elastic basis and rigidly supported alongside the margins. A general second-order obstacle problem can be given as:

$$u^{(4)}(x) = \begin{cases} g(x), a \leq x < c \\ f(x)u(x) + g(x) + r, c \leq x < d \\ g(x), d \leq x < b \end{cases} \tag{4.1}$$

with the associated BCs given as $u(a) = u(b) = \alpha_f$, $u^{(1)}(a) = u^{(1)}(b) = \beta_f$ and $u^{(2)}(a) = u^{(2)}(b) = \gamma_f$. It is given that the operators $u^{(1)}(x), u^{(2)}(x)$ and $u^{(3)}(x)$ remain constant on $c$ and $d$.

## 2. Formulation of Obstacle Problem

For better understanding, we will start with some 2$^{nd}$ order obstacle BVPs for finding the solution $u(x)$ given that

$$\left.\begin{array}{lll} -u''(x) \geq f(x) & \text{on} & \Omega = [0,1] \\ u(x) \geq \Psi(x) & \text{on} & \Omega = [0,1] \\ [u''(x) + f(x)][u(x) - \Psi(x)] = 0 & \text{on} & \Omega = [0,1] \\ u(0) = 0, u'(0) = 0, u'(1) = 0 \end{array}\right\} \quad (5)$$

Where $\Psi(x)$ and $f(x)$ represents the elastic obstacle function and continuous function. The obstacle equations (5) define the static form of the adjustable string, dragged eventually and stretched upon an flexible impediment $\Psi(x)$. To look at the problem in the light of VIT, we have to proceed in the following way. We construct the set $\mathfrak{K}$ as

$$\mathfrak{K} = \{v : v \in H_0^2(\Omega) \wedge v \geq \Psi \text{ on } \Omega\}$$

which represents a closed convex subset in $H_0^2(\Omega)$, where $H_0^2(\Omega)$ is Sobolove Space and Hilbert Space as well [8, 9]. The associated power operational $I[v]$ by impediment equation (5) through Tonti approach[10] like:

$$\begin{aligned} I[v] &= \int_0^1 (-v''(x) - 2f(x))v(x)dx \,\forall v \in \mathfrak{K} \\ &= -[v'(x)v(x)]|_0^1 + \int_0^1 v'(x)v'(x)dx - 2\int_0^1 f(x)v(x)dx \\ &= -[v(1)v'(1)] + [v(0)v'(0)] + \int_0^1 v'(x)v'(x)dx - 2\int_0^1 f(x)v(x)dx \\ &= \int_0^1 v'(x)v'(x)dx - 2\int_0^1 f(x)v(x)dx \\ &= \langle Tv, v \rangle - 2\langle f, v \rangle, \end{aligned}$$

where

$$\langle Tv, v \rangle = \int_0^1 v'(x)v'(x)dx \quad (6)$$

and

$$\langle f, v \rangle = \int_0^1 f(x)v(x)dx. \quad (7)$$

We can finalize that the manipulator charactorized through (6) is symmetric, positive and linear. Each requirments for [7 Theorem 3.1, 8, 9 Chapter IV Section 2] are fulfilled. Therefore in that case a minimal $u(x)$ for operational $I[v]$ specified above, on the convex and closed set $\mathfrak{K}$ of $H_0^2(\Omega)$ can be described through a changing inequity of the type given below in (8) as the notion of Stampacchia and Lewy [7], equation(5) may be expressed like:

$$\left.\begin{array}{l} -u''(x) + \mu\{(u(x) - \Psi(x))\}(u(x) - \Psi(x)) = f, 0 < x < 1 \\ u(0) = 0, u'(0) = 0, u'(1) = 0 \end{array}\right\} \quad (8)$$

Where

$$\mu(t) = \begin{cases} -1, t \geq 0 \\ 0, t < 0 \end{cases} \quad (9)$$

represents the penalty function, and an obstacle problem $\Psi(x)$ is defined as:

$$\Psi(x) = \begin{cases} -1 & for\ 0 \leq x < \frac{1}{4} \\ 1 & for\ \frac{1}{4} \leq x < \frac{3}{4} \\ -1 & for\ \frac{3}{4} \leq x \leq 1 \end{cases} \quad (10)$$

From equations (8), (9), and (10), we obtained the following system of 2$^{nd}$ order differential equations

$$u''(x) = \begin{cases} f & for\ 0 \leq x < \frac{1}{4} \\ u(x) + f - 1 & for\ \frac{1}{4} \leq x < \frac{3}{4} \\ f & for\ \frac{3}{4} \leq x \leq 1 \end{cases} \quad (11)$$

with the BCs

$$u(0) = 0,\ u'(0) = 0,\ u'(1) = 0$$

and conditions that the function $u(x)$ and $u'(x)$ are continuous at $x = \frac{1}{4}$ and $x = \frac{3}{4}$.

### 2.1. Cantilinear Bridge Problem

The motivation for this problem is due to the bridge problems in civil engineering. Saint John River build in New Brunswick, the Niagara River between the state of New York and Ontario, the Fraser River in British Columbia, and the great Forth Bridge in Scotland are a few examples of cantilever bridges in the world. Howrah Bridge is a cantilever bridge with a suspended span. The Hooghly River in West Bengal, India is also an example of a cantilever bridge which is the largest bridge in the world.

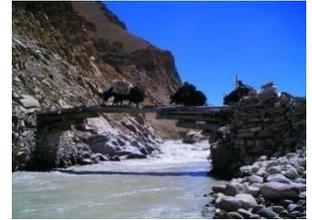

The formulation of this bridge is the invention of great minds of civil engineers who use the beauty of mathematics to formalize it as a mathematical problem [23]. To improve the system, it was essential to model the problem in such a way that we can easily exercise the practical situation on simulation to avoid and control unplanned situations. To perk up the accuracy of the model, durability, and sustainability of the structure and convergence of the multi-support bridge problem, the mathematicians design it as a linear obstacle problem using the concept of variational inequalities.

Now in order to enhance precision and convergency for the previously described equations modern technique has now been used for cable-connected also uninterrupted bridge multi function linearly impediment equations.

A modern technique has now been used for cable-connected also uninterrupted bridge multi function linearly impediment equations in order to enhance precision and convergency for the previously described equations. In comparison with ordinary and numeric techniques, a creative technique provides best also accurate solutions. The constant support bridge equations can effectively illuminate by means of inventive strategy, that is decaying method and semi-analytic approach.The deterioration inventive strategy comes about are found to focalize exceptionally rpidly and are more near to the precise approximation.There are numerous deterrents in multi-support equations such as cable-connected and constant support bridges; easily one can intepret the equation succesfully utilizing creative strategies.Illustrated comes about shortages that the creative strategy can utilize any sort of ordinary equation for which the approximated calculations are inconvenient.

The ever-changing disparities give common, novel and inventive setups in order to the detailing of the multi-support,partial,intercourse,economical and also maximization. In an inovative detailing of such equations, the area of the intercourse region open boundary gets to be an inborn portion of the result and no uncommon procedures are required to get it. In case the impediment is recognized, at that point the innovative imbalances may be classified by a framework of differential conditions by utilizing the retribution work strategy described by Stampacchia and Lewy [3]. The major calculational benefit of this method is it is straightforward pertinence for understanding the framework of differential conditions. During latest years, Noor and Al-Said [4], Noor and Khalifa [5], Al-Said and Noor [6], and Al-Said et al. [7] have utilized allocation, limited contrast, as well as spline strategies for tackling these kind of second order framework of differential problems related with multi-support as well as partial equations. Within the latest paper, the current approach is utilized of getting results of the nth-order boundary value approach of this kind.

In the next section, we will find the exact solution to the obstacle problems and discuss the difficulties in finding the exact solution.

### 3. Exact Solution of Obstacle Problems

In this section, we will find the different obstacle problems for different values of the continuous function $g(x)$, $f(x)$, $a, b, r, c, d$, and the corresponding boundary conditions. We will find exact solutions to the designed Obstacle Problems. This section will also view the graph of the exact solution of the considered Obstacle Problems. This section will also discover the estimate result of the considered Obstacle equations by considering OHAM and compare them with the exact solution.

### 3.1. Numerical Problems

**Example 3.1.1:** Consider the 2$^{nd}$ order obstacle problem given in section 1 in equation (3) represented as follows

$$u^{(2)}(x) = \begin{cases} g(x), d \leq x < b \\ f(x)u(x) + g(x) + r, c \leq x < d \\ g(x), a \leq x < c \end{cases}$$

with $g(x) = 0$, $f(x) = 1$, $r = -1$, $a = -1, c = -\frac{1}{2}, d = \frac{1}{2}$ and $b = 1$. The BCs given are $u(1) = u(-1) = 0$.

To find the precise result to the above example, we can proceed in the following way:

First, we consider the differential equations

$u_0^{(2)}(x) = 0, u_0(-1) = 0$ and solve with respect to $x$ for $u_0(x)$

$u_1^{(2)}(x) = u_1(x) - 1$, and solve with respect to $x$ for $u_1(x)$

$u_2^{(2)}(x) = 0, u_2(1) = 0$ and solve with respect to $x$ for $u_2(x)$

By solving each differential equation by the method of undetermined coefficients we will get

$u_0(x) = (1 + x)a_1$               (p1.1)

$u_1(x) = 1 + a_2 e^x + a_3 e^{-x}$      (p1.2)

$u_2(x) = a_4(1 - x)$              (p1.3)

In this system of solutions, we have to find the constants $a_1, a_2, a_3$ and $a_4$ using the condition the functions $u^{(1)}(x)$ is continuous on $c$ and $d$. To compute the values of $a_1, a_2, a_3$ and $a_4$, Gauss elimination method is used as:

$$u_0'\left(-\frac{1}{2}\right) - u_1'\left(-\frac{1}{2}\right) = 0$$

$$u_1'\left(\frac{1}{2}\right) - u_2'\left(\frac{1}{2}\right) = 0 \tag{p1.4}$$

After finding the values, accurate results are as follows:

$$u(x) = \begin{cases} \frac{2(-1+e)(1+x)}{1+3e}, & a \leq x < c \\ 1 - \frac{4\sqrt{e}\cosh[x]}{1+3e}, & c \leq x < d \\ -\frac{2(-1+e)(-1+x)}{1+3e}, & d \leq x < b \end{cases} \tag{p1.5}$$

The graphical view of the precise approximation emerges as

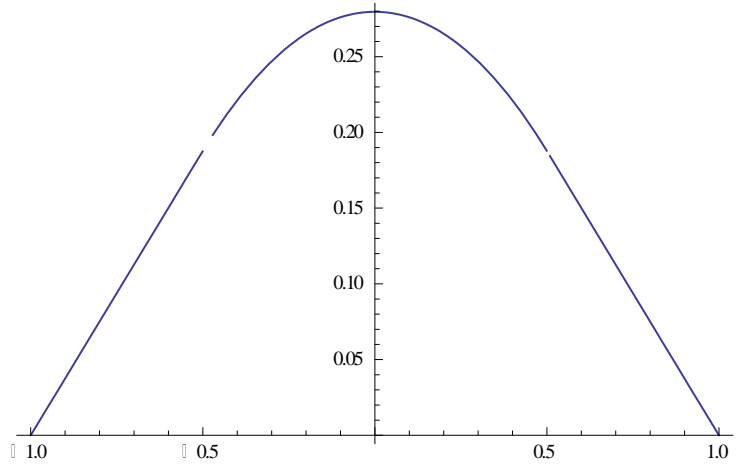

**Example 3.1.2:** Consider the 2nd order obstacle problem given in section 1 in equation (3) represented as follows

$$u^{(2)}(x) = \begin{cases} g(x), & a \leq x < c \\ f(x)u(x) + g(x) + r, & c \leq x < d \\ g(x), & d \leq x < b \end{cases}$$

with $g(x) = x$, $f(x) = 1$, $r = -1$, $a = 0$, $c = \frac{1}{4}$, $d = \frac{3}{4}$ and $b = 1$. The BCs given are $u(0) = u(1) = 0$.

Working as in example 3.1.1, the exact solution is given as

$$u(x) = \begin{cases} -\frac{(-2049+80\sqrt{e}+545e)x}{96(-9+25e)} + \frac{x^3}{6}, & a \leq x < c \\ 1 + \frac{e^{-\frac{1}{4}x}\left(15e - 965e^{3/2} + 579e^{2x} - 25e^{\frac{1}{2}+2x}\right)}{48(-9+25e)} - x, & c \leq x < d \\ \frac{1}{6}(-1+x)\left(\frac{933+3088\sqrt{e}-2725e}{-144+400e} + x + x^2\right), & d \leq x < b \end{cases} \tag{p2.1}$$

The graphical view of the exact solution is given as

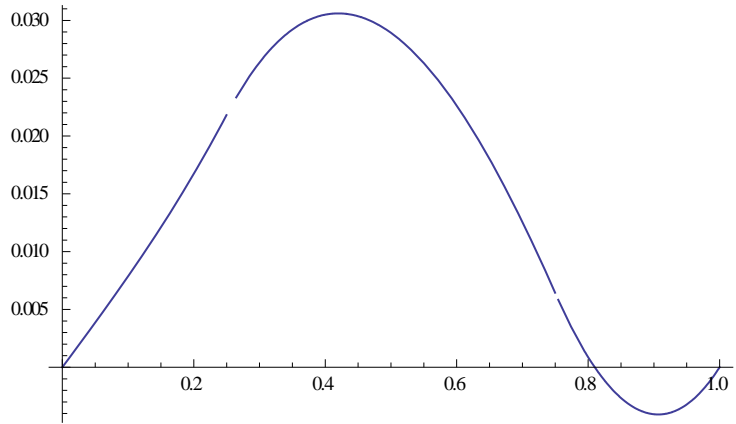

**Example 3.1.3:** Consider the 2$^{nd}$ order obstacle problem given in section 1 in equation (3) represented as follows

$$u^{(2)}(x) = \begin{cases} g(x), a \leq x < c \\ f(x)u(x) + g(x) + r, c \leq x < d \\ g(x), d \leq x < b \end{cases}$$

with $g(x) = u'(x) - 2$, $f(x) = 1$, $r = -1$, $a = 0$, $c = \frac{1}{4}$, $d = \frac{3}{4}$ and b= 1. And boundary constraints given are $u(1) = u(0) = 0$.

Working as in example 3.1.1, accurate results are as follows:

$$u(x) = \begin{cases} 2x + (-1 + e^x)a_1, a \leq x < c \\ 3 + e^{-\frac{1}{2}(-1+\sqrt{5})x}\left(f + ge^{\sqrt{5}x}\right), c \leq x < d \\ -2 + 2x + (-e + e^x)a_2, d \leq x < b \end{cases} \quad \text{(p3.1)}$$

where

$$a_1 = \frac{-11\sqrt{5} + e^{1/4}\left(4\sqrt{5} + \sqrt{5}(9 - 4e^{1/4})\cosh\left[\frac{\sqrt{5}}{4}\right] + (-19 + 6e^{1/4})\sinh\left[\frac{\sqrt{5}}{4}\right]\right)}{2e^{1/4}\left(\sqrt{5}(-1 + \sqrt{e})\cosh\left[\frac{\sqrt{5}}{4}\right] + (3 - 4e^{1/4} + 3\sqrt{e})\sinh\left[\frac{\sqrt{5}}{4}\right]\right)}$$

$f$

$$= \frac{\left(e^{-\frac{5}{8}+\frac{3\sqrt{5}}{8}}\left(-(-11 + 4e^{1/4})(5 + \sqrt{5} + (-5 + \sqrt{5})e^{1/4})e^{1/4} - (-4 + 9e^{1/4})(-5 + \sqrt{5} + (5 + \sqrt{5})e^{1/4})e^{\frac{1}{4}(}\right.\right.}{\left(2\sqrt{5}\left(-3 - \sqrt{5} + 4e^{1/4} - 3\sqrt{e} + (3 - \sqrt{5} - 4e^{1/4} + (3 + \sqrt{5})\sqrt{e})e^{\frac{\sqrt{5}}{2}} + \sqrt{5e}\right)\right)}$$

$$g = \frac{\left(e^{\frac{1}{8}(-5-\sqrt{5})}\begin{pmatrix}-4(1+\sqrt{5})\sqrt{e}+(5+13\sqrt{5})e^{\frac{3}{4}}-9(-1+\sqrt{5})e+\\ \left(-11+4e^{\frac{1}{4}}\right)\left(1-\sqrt{5}+(1+\sqrt{5})e^{\frac{1}{4}}\right)e^{\frac{1}{4}(1+\sqrt{5})}\end{pmatrix}\right)}{\left(2\left(-3-\sqrt{5}+4e^{\frac{1}{4}}-3\sqrt{e}+\left(3-\sqrt{5}-4e^{\frac{1}{4}}+(3+\sqrt{5})\sqrt{e}\right)e^{\frac{\sqrt{5}}{2}}+\sqrt{5}e\right)\right)}$$

$$a_2 = \frac{\left(-10+4\sqrt{5}+(25-11\sqrt{5})e^{1/4}+(10+4\sqrt{5}-(25+11\sqrt{5})e^{1/4})e^{\frac{\sqrt{5}}{2}}+2\sqrt{5}(-4+9e^{1/4})e^{\frac{1}{4}(1+\sqrt{5})}\right)}{\left(-2(3+\sqrt{5})e^{3/4}+8e+2(-3+\sqrt{5})e^{5/4}+2(3-\sqrt{5}-4e^{1/4}+(3+\sqrt{5})\sqrt{e})e^{\frac{3}{4}+\frac{\sqrt{5}}{2}}\right)}$$

The graphical view of the exact solution is given as

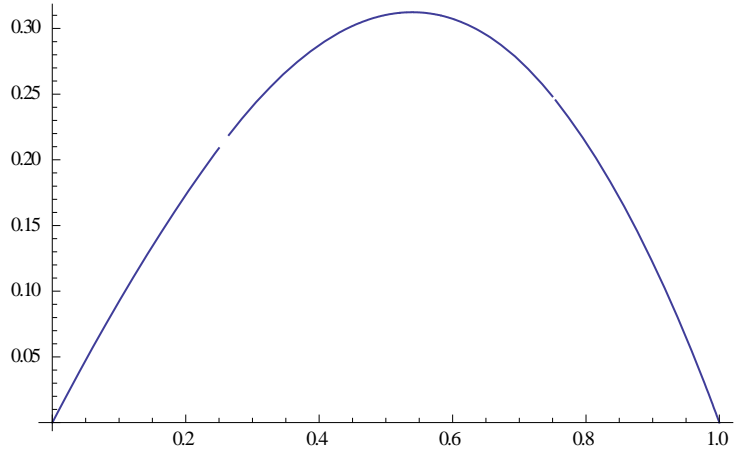

**Example 3.1.4:** Consider the 2nd order obstacle problem given in section 1 in equation (3) represented as follows

$$u^{(2)}(x) = \begin{cases} g(x), \eth \le x < b \\ f(x)u(x) + g(x) + r, c \le x < \eth \\ g(x), a \le x < c \end{cases}$$

with $g(x) = 0, f(x) = 1, r = -1, a = 0, c = \frac{\pi}{4}, \eth = \frac{3\pi}{4}$ and $b = \pi$. The BCs given are $u(0) = u(1) = 0$.

Working as in example 3.1.1, the exact solution is:

$$u(x) = \begin{cases} a_1 x, a \le x < c \\ 1 + e^x a_2 + e^{-x} a_3, c \le x < \eth \\ \frac{a_4(\pi - x)}{\pi}, \eth \le x < b \end{cases} \quad (p4.1)$$

Where

$$a_1 = \frac{4}{\pi + 4\coth\left[\frac{\pi}{4}\right]}$$

$$a_2 = -\frac{4e^{-\pi/4}}{4 - \pi + e^{\pi/2}(4 + \pi)}$$

$$a_3 = -\frac{4e^{3\pi/4}}{4 - \pi + e^{\pi/2}(4 + \pi)}$$

$$a_4 = \frac{4\pi}{\pi + 4\coth\left[\frac{\pi}{4}\right]}$$

The graphical view of the accurat results are as follows:

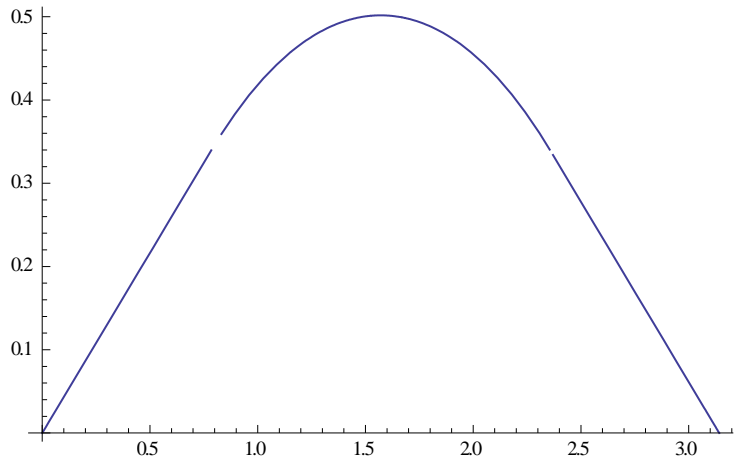

**Example 3.1.5:** Suppose the 2$^{nd}$ order obstacle problem given in section 1 in equation (3) represented as follows

$$-u^{(2)}(x) = \begin{cases} g(x), \mathfrak{d} \leq x < \mathfrak{b} \\ \mathfrak{f}(x)u(x) + g(x) + \mathfrak{r}, \mathfrak{c} \leq x < \mathfrak{d} \\ g(x), \mathfrak{a} \leq x < \mathfrak{c} \end{cases}$$

with $g(x) = u + 1$, $\mathfrak{f}(x) = 1$, $\mathfrak{r} = -1$, $\mathfrak{a} = 0$, $\mathfrak{c} = \frac{\pi}{4}$, $\mathfrak{d} = \frac{3\pi}{4}$ and $\mathfrak{b} = \pi$. The boundary-terms given are $u(0) = u(1) = 0$.

Working as in example 3.1.1, the exact solution is given as

$$u(x) = \begin{cases} -1 + \cos[x] + \left(-1 + \sqrt{2} + \sqrt{2}k\right)\sin[x], \mathfrak{a} \leq x < \mathfrak{c} \\ \dfrac{3k\pi - p\pi - 4kx + 4px}{2\pi}, \mathfrak{c} \leq x < \mathfrak{d} \\ -1 - \cos[x] + \left(-1 + \sqrt{2} + \sqrt{2}p\right)\sin[x], \mathfrak{d} \leq x < \mathfrak{b} \end{cases}$$

Where

$$k = -1 + \sqrt{2}$$

$$p = -1 + \sqrt{2}$$

The graphical view of the exact solution is given as

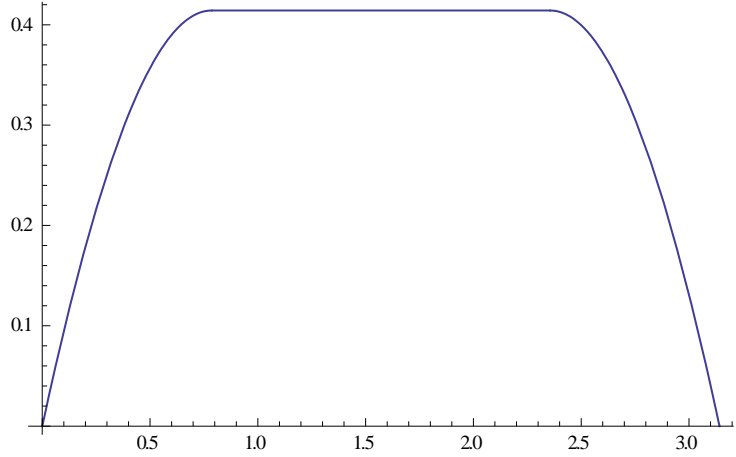

**Example 3.1.6:** Assume the 3$^{rd}$ order obstacle problem given in section 1 in equation (3) represented as follows

$$u^{(3)}(x) = \begin{cases} g(x), \mathfrak{d} \leq x < b \\ f(x)u(x) + g(x) + r, c \leq x < \mathfrak{d} \\ g(x), a \leq x < c \end{cases}$$

with $g(x) = 0, f(x) = 1, r = -1, a = 0, c = \frac{1}{4}, \mathfrak{d} = \frac{3}{4}$ and b= 1. And boundary constraints given are $u(1) = u(0) = 0,\ u'(\frac{3}{4}) = u'(\frac{1}{4}) = 0$.

To compute the exact solution to problem 5 we can proceed in the following way:

First, we consider the differential equations

$u_0^{(3)}(x) = 0, u_0(1) = 0$ and solve with respect to $x$ for $u_0(x)$

$u_1^{(3)}(x) = u_1(x) - 1$, and solve with respect to $x$ for $u_1(x)$

$u_2^{(3)}(x) = 0, u_2(1) = 0$ and solve with respect to $x$ for $u_2(x)$

By solving each differential equation by the method of undetermined coefficients we will get

$u_0(x) = x(a + bx)$ \hfill (p1.1)

$u_1(x) = 1 + e^x t + e^{-x/2}\left(p \cos\left[\frac{\sqrt{3}x}{2}\right] + q \sin\left[\frac{\sqrt{3}x}{2}\right]\right)$ \hfill (p1.2)

$u_2(x) = (1-x)(dx + c(1+x))$ \hfill (p1.3)

In this system of solutions, we have to find the constants $a,b,c,d,t,q,p$ using the condition of the functions $u^{(1)}(x)$ is continuous on c and d. To calculate $a, b, c, d, t, q$, and $p$, Gauss elimination method is used as:

$$u_0'[1/4] - u_1'[1/4] == 0, u_1'[3/4] - u_2'[3/4] == 0$$
$$u_0''[1/4] - u_1''[1/4] == 0, u_1''[3/4] - u_2''[3/4] == 0 \hfill (p1.4)$$

After finding the values, accurate results are as follows:

$$u(x) = \begin{cases} x(a+bx), a \leq x < c \\ 1 + e^x t + e^{-x/2}\left(p\cos\left[\frac{\sqrt{3}x}{2}\right] + q\sin\left[\frac{\sqrt{3}x}{2}\right]\right), c \leq x < d \\ (1-x)(dx + c(1+x)), d \leq x < b \end{cases} \qquad (p3.1)$$

Where if $c = 1$ we can have

$$a = \frac{\begin{pmatrix} 888 + 41(24-41c)e + (-984 + 1297c - 888\sqrt{e})e^{1/4}\cos\left[\frac{\sqrt{3}}{4}\right] + \\ \sqrt{3}(1048 + 527c - 1208\sqrt{e})e^{1/4}\sin\left[\frac{\sqrt{3}}{4}\right] \end{pmatrix}}{\left(-925 + e^{3/4}\left(1309\cos\left[\frac{\sqrt{3}}{4}\right] - 1059\sqrt{3}\sin\left[\frac{\sqrt{3}}{4}\right]\right)\right)}$$

$$b = \frac{\left(4\left(148 + 5(-24 + 41c)e + (-264 - 13c + 236\sqrt{e})e^{\frac{1}{4}}\cos\left[\frac{\sqrt{3}}{4}\right] + \sqrt{3}(8 + 117c + 44\sqrt{e})e^{\frac{1}{4}}\sin\left[\frac{\sqrt{3}}{4}\right]\right)\right)}{\left(-925 + e^{\frac{3}{4}}\left(1309\cos\left[\frac{\sqrt{3}}{4}\right] - 1059\sqrt{3}\sin\left[\frac{\sqrt{3}}{4}\right]\right)\right)}$$

$$q = \frac{\left(e^{3/8}(B_1 - B_2 + B_3)\right)}{\left(-2775 + e^{3/4}\left(3927\cos\left[\frac{\sqrt{3}}{4}\right] - 3177\sqrt{3}\sin\left[\frac{\sqrt{3}}{4}\right]\right)\right)}$$

where

$$B_1 = \sqrt{3}(600 + 1375c - 2304\sqrt{e})\cos\left[\frac{3\sqrt{3}}{8}\right]$$

$$B_2 = (-24 + 41c)e^{3/4}\left(71\sqrt{3}\cos\left[\frac{\sqrt{3}}{8}\right] + 27\sin\left[\frac{\sqrt{3}}{8}\right]\right)$$

$$B_3 = 75(40 - 7c)\sin\left[\frac{3\sqrt{3}}{8}\right] - 4416\sqrt{e}\sin\left[\frac{3\sqrt{3}}{8}\right]$$

$$d = \frac{400 + 1225c - 1536\sqrt{e} + (1136 - 2377c)e^{\frac{3}{4}}\cos\left[\frac{\sqrt{3}}{4}\right] + \sqrt{3}(-144 + 599c)e^{\frac{3}{4}}\sin\left[\frac{\sqrt{3}}{4}\right]}{-925 + e^{\frac{3}{4}}\left(1309\cos\left[\frac{\sqrt{3}}{4}\right] - 1059\sqrt{3}\sin\left[\frac{\sqrt{3}}{4}\right]\right)}$$

$$t = \frac{3552 + 2e^{1/4}\left(24(-82 + 31c)\cos\left[\frac{\sqrt{3}}{4}\right] + \sqrt{3}(696 + 929c)\sin\left[\frac{\sqrt{3}}{4}\right]\right)}{-2775e^{1/4} + 3927e\cos\left[\frac{\sqrt{3}}{4}\right] - 3177\sqrt{3}e\sin\left[\frac{\sqrt{3}}{4}\right]}$$

$$p = \frac{\left(e^{3/8}(A1 + A2 + A3)\right)}{\left(-2775 + e^{3/4}\left(3927\cos\left[\frac{\sqrt{3}}{4}\right] - 3177\sqrt{3}\sin\left[\frac{\sqrt{3}}{4}\right]\right)\right)}$$

where

$$A1 = 75(40 - 7c)\cos\left[\frac{3\sqrt{3}}{8}\right] - 4416\sqrt{e}\cos\left[\frac{3\sqrt{3}}{8}\right]$$

$$A2 = (-24 + 41c)e^{3/4}\left(-27\cos\left[\frac{\sqrt{3}}{8}\right] + 71\sqrt{3}\sin\left[\frac{\sqrt{3}}{8}\right]\right)$$

$$A3 = \sqrt{3}(-600 - 1375c + 2304\sqrt{e})\sin\left[\frac{3\sqrt{3}}{8}\right]$$

The graphical view of the exact solution is given as

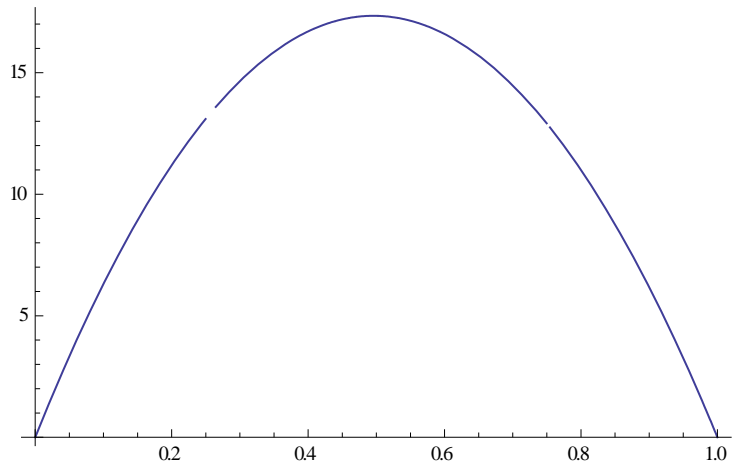

**Example 3.1.7:** Suppose the 3$^{rd}$ order obstacle problem given in section 1 in equation (3) represented as follows

$$u^{(3)}(x) = \begin{cases} g(x), d \leq x < b \\ f(x)u(x) + g(x) + r, c \leq x < d \\ g(x), a \leq x < c \end{cases}$$

with $g(x) = 2, f(x) = 1, r = 1, a = 0, c = \frac{1}{4}, d = \frac{3}{4}$ and b= 1. The BCs given are $u(0) = u(1) = 0$, $u'(\frac{1}{4}) = u'(\frac{3}{4}) = 0$.

After finding the values, the exact solution is given as

$$u(x) = \begin{cases} \frac{1}{3}x(3a + 3bx + x^2), a \leq x < c \\ 1 + e^x t + e^{-x/2}\left(p\cos\left[\frac{\sqrt{3}x}{2}\right] + q\sin\left[\frac{\sqrt{3}x}{2}\right]\right), c \leq x < d \\ \frac{1}{3}(-1 + x)(-3c - 3(c + d)x + x^2), d \leq x < b \end{cases}$$

Where if $c = 1$ we can have

$$a = \frac{-(39627 + X_1 + X_2 + X_3)}{\left(48\left(925 + e^{3/4}\left(-1309\cos\left[\frac{\sqrt{3}}{4}\right] + 1059\sqrt{3}\sin\left[\frac{\sqrt{3}}{4}\right]\right)\right)\right)}$$

Where

$$X_1 = 41(871 - 1968c)e$$

$$X_2 = \left(-38399 + 62256c - 38475\sqrt{e}\right)e^{1/4}\cos\left[\frac{\sqrt{3}}{4}\right]$$

$$X_3 = \sqrt{3}\left(54255 + 25296c - 60859\sqrt{e}\right)e^{1/4}\sin\left[\frac{\sqrt{3}}{4}\right]$$

$$b = \frac{-(9842 - Y_1 + Y_2 + Y_3)}{\left(12\left(925 + e^{3/4}\left(-1309\cos\left[\frac{\sqrt{3}}{4}\right] + 1059\sqrt{3}\sin\left[\frac{\sqrt{3}}{4}\right]\right)\right)\right)}$$

Where

$$Y_1 = (4355 + 9840c)e$$

$$Y_2 = \left(-12829 - 624c + 7342\sqrt{e}\right)e^{1/4}\cos\left[\frac{\sqrt{3}}{4}\right]$$

$$Y_3 = \sqrt{3}\left(1205 + 5616c + 5278\sqrt{e}\right)e^{1/4}\sin\left[\frac{\sqrt{3}}{4}\right]$$

$$t = \frac{84804 + e^{1/4}\left(36(-2493 + 992c)\cos\left[\frac{\sqrt{3}}{4}\right] + \sqrt{3}(40085 + 44592c)\sin\left[\frac{\sqrt{3}}{4}\right]\right)}{72\left(-925e^{1/4} + 1309e\cos\left[\frac{\sqrt{3}}{4}\right] - 1059\sqrt{3}e\sin\left[\frac{\sqrt{3}}{4}\right]\right)}$$

$$p = \frac{\left(e^{3/8}(W_1 + W_2 + W_3)\right)}{\left(144\left(925 + e^{3/4}\left(-1309\cos\left[\frac{\sqrt{3}}{4}\right] + 1059\sqrt{3}\sin\left[\frac{\sqrt{3}}{4}\right]\right)\right)\right)}$$

Where

$$W_1 = 225(627 - 112c)\cos\left[\frac{3\sqrt{3}}{8}\right]$$

$$W_2 = -210864\sqrt{e}\cos\left[\frac{3\sqrt{3}}{8}\right]$$

$$W_3 = (-871 + 1968c)e^{3/4}\left(-27\cos\left[\frac{\sqrt{3}}{8}\right] + 71\sqrt{3}\sin\left[\frac{\sqrt{3}}{8}\right]\right)$$

$$q = \frac{\left(e^{3/8}(Z_1 + Z_2 + Z_3 + Z_4)\right)}{\left(48\left(-925\sqrt{3} + e^{3/4}\left(1309\sqrt{3}\cos\left[\frac{\sqrt{3}}{4}\right] - 3177\sin\left[\frac{\sqrt{3}}{4}\right]\right)\right)\right)}$$

$$Z_1 = 25(1543 + 2640c)\cos\left[\frac{3\sqrt{3}}{8}\right]$$

$$Z_2 = -110016\sqrt{e}\cos\left[\frac{3\sqrt{3}}{8}\right]$$

$$Z_3 = (-871 + 1968c)e^{3/4}\left(71\cos\left[\frac{\sqrt{3}}{8}\right] + 9\sqrt{3}\sin\left[\frac{\sqrt{3}}{8}\right]\right)$$

$$Z_4 = \left[\sqrt{3}(47025 - 8400c - 70288\sqrt{e})\sin\left[\frac{3\sqrt{3}}{8}\right]\right]$$

$$d = \frac{\left(\sqrt{3}(725 + 14700c - 18336\sqrt{e}) + \sqrt{3}(18379 - 28524c)e^{3/4}\cos\left[\frac{\sqrt{3}}{4}\right] + 3(-7837 + 7188c)e^{3/4}\sin\left[\frac{\sqrt{3}}{4}\right]\right)}{\left(12\left(-925\sqrt{3} + e^{3/4}\left(1309\sqrt{3}\cos\left[\frac{\sqrt{3}}{4}\right] - 3177\sin\left[\frac{\sqrt{3}}{4}\right]\right)\right)\right)}$$

The graphical view of the exact solution is given as

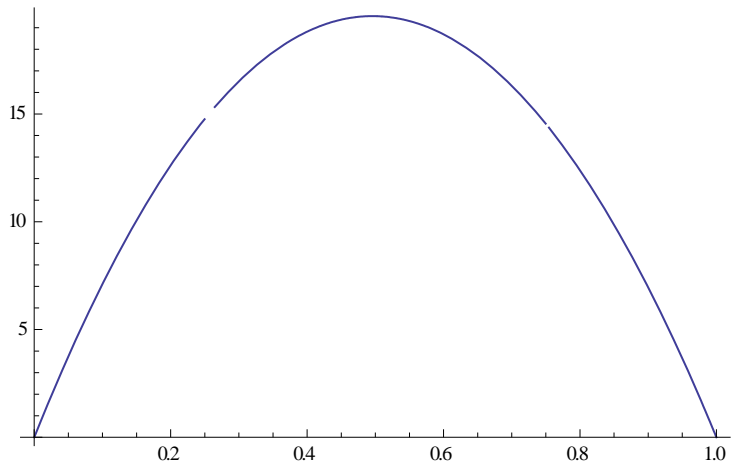

**Example 3.1.8:** Assume the 3<sup>rd</sup> order obstacle problem given in section 1 in equation (3) represented as follows

$$u^{(3)}(x) = \begin{cases} g(x), a \leq x < c \\ f(x)u(x) + g(x) + r, c \leq x < d \\ g(x), d \leq x < b \end{cases}$$

With $g(x) = x$, $f(x) = 1$, $r = -1$, $a = 0$, $c = \frac{1}{4}$, $d = \frac{3}{4}$ and $b = 1$. And boundary constraints given are $u(1) = u(0) = 0$, $u'(\frac{3}{4}) = u'(\frac{1}{4}) = 0$.

After finding the values, the exact solution is given as

$$u(x) = \begin{cases} a_1 x + a_2 x^2 + \dfrac{x^4}{24}, a \leq x < c \\ 1 + e^x a_3 - x + e^{-x/2}\left(a_4 \cos\left[\dfrac{\sqrt{3}x}{2}\right] + a_5 \sin\left[\dfrac{\sqrt{3}x}{2}\right]\right), c \leq x < d \\ \dfrac{1}{24}(-1 + x)(-24 a_6 x + (1 + x)(-24 a_7 + x^2)), d \leq x < b \end{cases}$$

where if $a_7 = 1$ we can have

$$a_1 = \frac{-(1388351 + L_1 + L_2 - L_3)}{\left(768\left(925 + e^{3/4}\left(-1309 \cos\left[\dfrac{\sqrt{3}}{4}\right] + 1059\sqrt{3} \sin\left[\dfrac{\sqrt{3}}{4}\right]\right)\right)\right)}$$

Where

$$L_1 = 41(30215 - 31488c)e$$

$$L_2 = (-956191 + 996096c - 1681727\sqrt{e})e^{1/4} \cos\left[\dfrac{\sqrt{3}}{4}\right]$$

$$L_3 = \left[3\sqrt{3}\left(-134912c + 2627(49 + 15\sqrt{e})\right)e^{1/4} \sin\left[\dfrac{\sqrt{3}}{4}\right]\right]$$

$$a_2 = \frac{-(232767 + M_1 + M_2 + M_3)}{\left(384\left(925 + e^{3/4}\left(-1309 \cos\left[\dfrac{\sqrt{3}}{4}\right] + 1059\sqrt{3} \sin\left[\dfrac{\sqrt{3}}{4}\right]\right)\right)\right)}$$

Where

$$M_1 = (-302150 + 314880c)e$$

$$M_2 = (18278 - 19968c + 354465\sqrt{e})e^{1/4} \cos\left[\dfrac{\sqrt{3}}{4}\right]$$

$$M_3 = 13\sqrt{3}\left(-13246 + 13824c + 5685\sqrt{e}\right)e^{1/4} \sin\left[\frac{\sqrt{3}}{4}\right]$$

$$a_3 = \frac{1363302 + e^{1/4}\left(18(-30599 + 31744c) \cos\left[\frac{\sqrt{3}}{4}\right] + \sqrt{3}(-682613 + 713472c) \sin\left[\frac{\sqrt{3}}{4}\right]\right)}{1152\left(-925e^{1/4} + 1309e \cos\left[\frac{\sqrt{3}}{4}\right] - 1059\sqrt{3}e \sin\left[\frac{\sqrt{3}}{4}\right]\right)}$$

$$a_4 = \frac{-\left(e^{3/8}(N_1 + N_2 + N_3 + N_4)\right)}{\left(2304\left(925 + e^{3/4}\left(-1309 \cos\left[\frac{\sqrt{3}}{4}\right] + 1059\sqrt{3} \sin\left[\frac{\sqrt{3}}{4}\right]\right)\right)\right)}$$

Where

$$N_1 = 225(1739 - 1792c) \cos\left[\frac{3\sqrt{3}}{8}\right]$$

$$N_2 = -3389832\sqrt{e} \cos\left[\frac{3\sqrt{3}}{8}\right]$$

$$N_3 = (-30215 + 31488c)e^{3/4}\left(-27 \cos\left[\frac{\sqrt{3}}{8}\right] + 71\sqrt{3} \sin\left[\frac{\sqrt{3}}{8}\right]\right)$$

$$N_4 = \sqrt{3}\left(1011025 - 1056000c + 1768608\sqrt{e}\right) \sin\left[\frac{3\sqrt{3}}{8}\right]$$

$$a_5 = \frac{\left(e^{3/8}(O_1 + O_2 + O_3 + O_4)\right)}{\left(768\left(-925\sqrt{3} + e^{3/4}\left(1309\sqrt{3} \cos\left[\frac{\sqrt{3}}{4}\right] - 3177 \sin\left[\frac{\sqrt{3}}{4}\right]\right)\right)\right)}$$

Where

$$O_1 = 25(-40441 + 42240c) \cos\left[\frac{3\sqrt{3}}{8}\right]$$

$$O_2 = -1768608\sqrt{e} \cos\left[\frac{3\sqrt{3}}{8}\right]$$

$$O_3 = -(-30215 + 31488c)e^{3/4}\left(71 \cos\left[\frac{\sqrt{3}}{8}\right] + 9\sqrt{3} \sin\left[\frac{\sqrt{3}}{8}\right]\right)$$

$$O_4 = -3\sqrt{3}\left(-43475 + 44800c + 376648\sqrt{e}\right) \sin\left[\frac{3\sqrt{3}}{8}\right]$$

$$a_6 = \frac{D_1 + D_2 + D_3}{\left(384\left(-925\sqrt{3} + e^{3/4}\left(1309\sqrt{3}\cos\left[\frac{\sqrt{3}}{4}\right] - 3177\sin\left[\frac{\sqrt{3}}{4}\right]\right)\right)\right)}$$

Where

$$D_1 = \sqrt{3}(-145225 + 470400c - 589536\sqrt{e})$$

$$D_2 = \sqrt{3}(443689 - 912768c)e^{3/4}\cos\left[\frac{\sqrt{3}}{4}\right]$$

$$D_3 = 3(128921 + 230016c)e^{3/4}\sin\left[\frac{\sqrt{3}}{4}\right]$$

The graphical view of the exact solution is given as

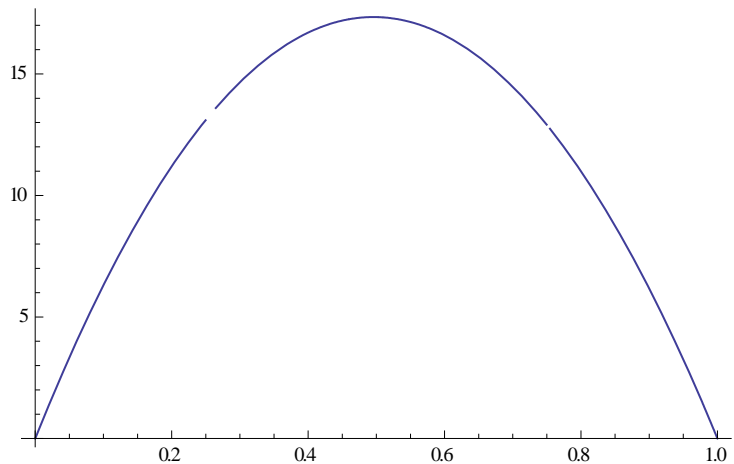

It is not always easy to find the exact solution of every choice of functions $f(x)$, $g(x)$, constant value $r$, $a, c, d$, and $b$ and under the boundary conditions given are $u(a) = u(b) = k_i$, $u'(c) = u'(d) = m_i$. When the functions involved are nonlinear then finding the exact solution is a very difficult job. This problem originates the formation of some numerical technique that can be used to solve at most exactly as we required within the situation.

## 4. Conclusion

A brief introduction about the obstacle problems is given and the problem is formulated. The general form is formulated up to the 3$^{rd}$ order. The exact solutions to the obstacle problems are calculated for several numerical problems in Mathematica. The graphical form of the solution is also presented for each problem. The graphical representation of results depicts the existence of symmetry among the results.

# References

1. A. Khan and T. Aziz(2003), "The numerical solution of third-order boundary-value problems using quintic splines," *Applied Mathematics and Computation, vol. 137, no. 2-3, pp. 253–260.*
2. A. K. Khalifa and M. A. Noor(1990), Quintic splines solutions of contact problems, *Math Comput Modell 13 51–58.*
3. A. Mohsen, M. El-Gamel(2008), On the Galerkin and collocation methods for two point boundary value problems using sinc bases, *Computers and Mathematics with Applications 56 930–941.*
4. A. Khan, I. Khan, T. Aziz(2006), Sextic spline solution of a singularly perturbed boundary-value problems, *Appl. Math. Comput. 181 432–439.*
5. A. Khan, T. Aziz(2003), Parametric cubic spline approach to the solution of a system of second-order boundary-value problems, *J. Optim. Theory Appl. 118 (1) 45–54.*
6. D. Kinderlehrer, G. Stampacchia(1980), An Introduction to Variational Inequalities and Their Applications, *Academic Press, London,.*
    Tonti, E (1984) Variational formulation for every nonlinear problem. *Int. J. Eng. Sci. 22: pp. 1343-1371.*
7. E. A. Al-Said and M. A. Noor (1995),, Computational methods for fourth-order obstacle boundary-value problems, *Commun Appl Nonlinear Anal 2 73–83.*
8. E.A.Al-SaidandM.A.Noor (1998), Numerical solution of a system of fourth order boundary value problems, *Int J Comput Math 70 347–355.*
9. E. A. Al-Said(2001), "Numerical solutions for system of third-order boundary value problems," *International Journal of Computer Mathematics, vol. 78, no. 1, pp. 111–121.*
10. E. A. Al-Said and M. A. Noor(2003), "Cubic splines method for a system of third-order boundary value problems*," Applied Mathematics and Computation, vol. 142, no. 2-3, pp. 195–204.*
11. F. Gao and C.-M. Chi(2006), "Solving third-order obstacle problems with quartic B-splines," *Applied Mathematics and Computation, vol. 180, no. 1, pp. 270–274.*
12. H. Lewy and G. Stampacchia(1969), On the regularity of the solution of the variational inequalities, *Commun Pure Appl Math 22, 153–188.*
13. H. Lewy, G. Stampacchia(1969), On the regularity of the solution of the variational inequalities, *Comm. Pure Appl. Math. 22 153– 188.*
14. H. N. Caglar, S. H. Caglar, and E. H. Twizell(1999), "The numerical solution of third-order boundary-value problems with fourth-degree B-spline functions," *International Journal of Computer Mathematics, vol. 71, no. 3, pp. 373–381,.*
15. I. A. Tirmizi, E. H. Twizell, and Siraj-Ul-Islam(2005), "A numerical method for third-order non-linear boundary-value problems in engineering*," International Journal of Computer Mathematics, vol. 82, no. 1, pp. 103–109.*
16. J. Cranck, (1984) Free and moving boundary-value problems, *Clarendon Press, Oxford, UK.*
17. J. Rashidinia, R. Jalilian et al(2007), Sextic spline method for the solution of a system of obstacle problems, *Applied Mathematics and Computation 190 1669–1674.*
18. K. Keller(1976), Numerical Solutions of Two-Point Boundary Value Problems, *SIAM, Philadelphia*
19. M. Aslam Noor(1988), General variational inequalities, *Appl. Math. Lett. 1 119–122.*
20. M.T. Heath(2002), Scientific Computing: An Introductory Survey, McGraw-Hill, New York.
21. M.A. Noor, A.K. Khalifa(1987), Cubic splines collocation methods for unilateral problems*, Int. J. Engng. Sci. 25 1527–1530.*
22. M.A. Noor, S.I.A. Tirmizi(1988), Finite difference techniques for solving obstacle problems, *Appl. Math. Lett. 1 267–271.*
23. M.A. Noor, K.L. Noor, Th. Rassias(1993), Some aspects of variational inequalities, *J. Comput. Appl. Math. 47 285–312.*